\documentclass[review]{elsarticle}

\usepackage{lineno,hyperref}
\usepackage{amsmath}

\journal{Journal of \LaTeX\ Templates}

\newtheorem{theorem}{Theorem}

\begin{document}

\begin{frontmatter}

\title{Local reconstruction from the spherical 
Radon transform in 3D}

\author{Rafik Aramyan}
\address{}
\fntext[myfootnote]{The work was supported by the Science Committee of RA, in the frames of the research project 21AG‐1A045}

\ead{rafikaramyan@yahoo.com}

\author[mymainaddress]{Institute of Mathematics of NAS RA, Yerevan, Armenia}


\address[mymainaddress]{24/5 Mar. Baghramyan av., Yerevan}

\begin{abstract}
In this article we study the spherical mean Radon transform in $\mathbf R^3$ with detectors centered on a plane. We use the consistency method suggested by the author of this article for the inversion of the transform in 3D. A new iterative inversion formula is presented.
This formula has the benefit of being local and is suitable for practical reconstructions. The inversion of the spherical mean Radon transform is required in mathematical models in thermo- and photo-acoustic tomography, radar imaging, and others.
\end{abstract}


\end{frontmatter}


\section{Introduction and formulation of the problem}\label{1}

\noindent Computer tomography has had a great impact on medical diagnostics.
The base of x-ray Tomography is the classical Radon transform that maps a function to its integrals over straight lines.
Researchers using various types of physical signals have developed new methods for Tomography.
Thermoacoustic tomography (TAT) is the most successful method (see the articles such as \cite{AgQu}-\cite{AmKu}, \cite{Amm}-\cite{An}, \cite{Gel},  \cite{KuKu}-\cite{Kr}, \cite{Nat}, \cite{No}, \cite{Pal}).
The concept of TAT is: a short-duration electromagnetic wave is sent through a biological object
to trigger a thermoacoustic
response in the tissue. Different biological cells reaction depends on the amount of energy that absorbed.
The thermoelastic expansion creates acoustic wave which can be measured by detectors placed outside the object (we assume that the sound speed is constant).
Thus, one effectively measures the integrals  of the energy absorption distribution function $f$
over all spheres centered at the detectors locations. A great diagnostic tool can be discovered if the energy absorption distribution function $f$ is found. To recover $f$ one needs to invert the so-called spherical Radon transform
of $f$  that integrates a function over spheres with the centers at detectors.

\noindent We denote by $\mathbf R^n$  ($n\geq2$) the Euclidean $n$ - dimensional space. Let ${\mathbf S^{n-1}}$ be the $n-1$ dimensional unit sphere in  $\mathbf R^n$ with the center at the origin  $O\in\mathbf R^n$. By $\mathcal{C}^\infty$ we denote the class of real valued functions for which all orders of derivatives are continuous and by $S(P,t)$ we denote the sphere of radius $t > 0$  centered at $P\in\mathbf R^n$.

\noindent We study the following mathematical problem. For a real valued function $f$ supported completely in a region $G\subset\mathbf R^n$, we are interested in recovering $f$ from the mean value $Mf(P,t)$ of $f$ over spheres $S(P,t)$ centered on some set
$L$; that is, given $Mf(P,t)$ for all $P\in L$ and  $t>0$, we wish to recover $f$.

\noindent TAT reconstruction arises the following problems. For which sets $ L$ is the data collected by detectors placed along $L$ sufficient for unique reconstruction of $f$ (so called a set of injectivity) and what are the inversion formulas?

\noindent  Agranovsky and Quinto in \cite{AgQu} have proved several significant uniqueness results
for the spherical Radon transform.

\noindent Obviously any line L (or a hyperplane in higher dimensions) is a non-uniqueness set.
On the other hand (see \cite{AgQu}), if
$f$ is supported completely on one side of a hyperplane $e$ (the standard situation in TAT), it
is uniquely recoverable from its spherical means centered on $e$.

\noindent Exact inversion formulas for the spherical Radon transform are currently known for different geometries of detectors (also, formulas are known for the problem with incomplete data): in 2D for lines, circles, ellipse, in hight dimensions for boundaries of special domains, including spheres, cylinders, ellipsoids and hyperplanes ( \cite{AgKuKu}-\cite{AmKu}, \cite{An},\cite{Ara19}, \cite{Be}-\cite{XW}).

\noindent In this paper for a real valued function $f$  defined in $\mathbf R^3$ and supported completely on one side of a plane $e$, we are interested in recovering $f$ from the mean value of $f$ over spheres centered on $e$.
The article suggests a new approach what is called a
consistency method for the inversion of the spherical Radon transform in 3D with detectors on a plane. By means of the method
a new iterative inversion formula is found which gives an algorithm to
recover an unknown smooth function supported in a region from its spherical means over spheres
centered on a plane outstand the region. Our reconstruction formula has the benefit of being local (see Theorem 2 below). Also, note that in our research the support region can be unbounded.

\noindent The consistency method, suggested by the author of the paper, first was applied in \cite{Ara10} (see also \cite{Ara09}, \cite{Ara11}, \cite{Ara21}) to invert generalized Radon transform on the sphere. The method was already used in \cite{Ara19} to invert the spherical mean Radon transform in 2D with detectors on a line. In \cite{Ara19} a new iterative formula to recover an unknown function was found which has the benefit of being local.

\noindent Now we consider the spherical Radon transform in $\mathbf R^3$ over spheres
centered on the plane $\{z=0\}$. For a continuous function $f$ supported in the compact region $G\subset\mathbf R^3$ located in the positive half space $\{z>0\}$ we have
\begin{equation}\label{1}
Mf(p,q,t)=\frac{1}{4\pi}\int_{\mathbf S^{2}}f((p,q,0)+t\omega)\,d\omega,
\end{equation}
for $(p,q,0)\in\{z=0\}; \, t\in[0,\infty)$. Here we consider  the restriction of  $f$  onto the sphere $S(p,q,t)$ with center $(p,q,0)\in \{z=0\}$ and radius $t > 0$, $\omega\in\mathbf S^{2}$ is the unit vector corresponds a point on $S(p,q,t)$ and we integrate with respect to the spherical Lebesgue measure on ${\mathbf S^{2}}$. The value
$Mf(p,q,t)$ is the average of $f$ over the sphere $S(p,q,t)$ with center $(p,q,0)\in \mathbf R^3$ and radius $t > 0$.

\noindent In this article we apply the consistency method: for every $(p.q,t)\in {\mathbf R^{2}}\times[0,\infty)$ we consider equation (1) as an
integral equation on the sphere $S(p,q,t)$. The general solution of the integral equation we write in terms of Fourier series expansion with unknown coefficients. Then we seek the unknown coefficients to find a family of consistent solutions. Thus the problem of recovering a real valued function $f$  from the mean value ${\cal M} f $ of $f$ over spheres centered on $z=0$ we reduce
to finding a consistent solutions of integral equations (1).

\noindent Now we describe the inversion formula. In $\mathbf R^3$ consider usual cartesian system of coordinates. $Mf(p,q,t)$ is the average of $f$ over the sphere with center $(p,q,0)\in \{z=0\}$ and radius $t \geq 0$.

\noindent We define a sequence of
{\it standard} polynomials $Q_{n,i}$ defined on the interval $[0,1]$, where $n,i$ are integers and $0\leq i\leq n$:
\begin{equation}\label{8}
Q_{n,i}(t)=\sum_{j=1}^{n+i}d_j(n,i)\,t^{2j},\,\,\,\,\,\,\,t\in[0,1]\end{equation}
with coefficients $d_j(n,i)$ for $1\leq j\leq n+i$. There are recurrent algorithms by means of which one can find the coefficients $d_j(n,i)$.
We call $Q_{n,i}$ standard polynomials because their construction does not depend on $f$.

\begin{theorem} Let $f\in\mathcal{C}^\infty$ be a real valued function  $f$ supported in a compact region
$G\subset\mathbf R^3$ located in half space $\{z>0\}$. For $(x,y,z)\in G$ we have
\begin{eqnarray}f(x,y,z)=\lim_{n\to\infty}2\times\qquad\qquad\label{9}\\
\left((2n^2+3n+1) Mf(x,y,z)+\sum_{i=0}^{n}\int_{0}^{z}z^{2i-1}Q_{n,i}(\frac{u}{z})\Delta^{i}(Mf(x,y,u))du\right)\nonumber\end{eqnarray}
here $Mf(x,y,u)$ is the average of $f$ over $S(x,y,u)$, $\Delta$  is the Laplace operator in two dimensions  with respect to Cartesian coordinates $x$ and $y$.
\end{theorem}

\noindent Theorem 1 suggests a practical algorithm to reconstruct $f$.

\noindent It follows from (\ref{9}) that for $(x,y,z)\in G$ the value $f(x,y,z)$ has a local description. Thus we have the following consequence of Theorem 1.

\begin{theorem} For $(x,y,z)\in G$ the value $f(x,y,z)$ depends on values $Mf$ in a neighborhood of $P=(x,y,0)$ and $0\leq t\leq z$.
\end{theorem}

\noindent Also, it follows from Theorem 1 that using (\ref{9}) one can recover an unknown smooth function supported in a region (which can be unbounded) from its spherical means over spheres
centered on a plane outstand the region.

\begin{figure}
\center
\includegraphics{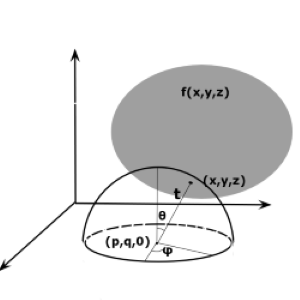}
\caption{An illustration to (1).}
\end{figure}

\section{Numerical Implementation and Computational example.}\label{S7}

The problem of reconstructing a function from spherical means is important for many
imaging and remote sensing applications (see, for example, \cite{FiRa}, \cite{KuKu}). These
applications require inversion of the  spherical
Radon transform. In some tomographic applications iterative reconstruction
algorithms are more common. In most cases, due to the presence of the derivative, the inversion formulas are sensitive to error in
the data $Mf$ (see \cite{Nat}, \cite{KuKu}).

\noindent In this paper, we suggest a new approach using the consistency method for the inversion of the spherical Radon transform in 3D. Our iterative formula to
recover an unknown smooth function $f$ supported in a compact region from
its spherical means $Mf$ (see Theorem 1) is different from the existing ones.
Our reconstruction formula has the benefit of being local.
Using iterative formula (\ref{9}) presented above for the step $n$ one can easily develop reconstruction algorithm. Computation of derivatives
in the formula (\ref{9}) can be easily done, for instance by using finite differences.
However, due to the presence of the derivatives of higher order, the inversion formulas are sensitive to error in the data $Mf$.

\noindent One can estimate the iteration speed using the following known result from the  theory of Fourier series expansion.

\textbf{Example 1}. Figure 2 provides a simple illustration of the iterative algorithm obtained by (\ref{9}) for a smooth function. We consider the following function

\begin{equation}\label{75}
f(x,y,z)=\begin{cases}
z^3y^2x \,\,\,\,\,\, \texttt{for}\,\,\, z\geq0\\
0\,\,\,\,\,\,  \texttt{for}\,\,\,
 z<0.\end{cases}
\end{equation}

\noindent Using Mathematica in Figure 2 we see the plot of the restriction of $f$ onto the plane $y=3$.

\begin{figure}
\center
\includegraphics{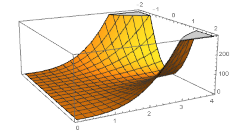}
\caption{the restriction of $f$ onto the plane $y=3$.}
\end{figure}

\noindent For average $Mf(x,y,u)$ of $f$ over the sphere with center
$(x,y,0)\in \{z=0\}$ and radius $u\geq 0$ we obtain
\begin{equation}\label{76}
Mf(x,y,u)=\frac{1}{8}x^2yu^3+\frac{1}{48}yu^5
\end{equation}
Also using recurrent relations for $n=2$ we obtain
\begin{equation}\label{77}\begin{cases}
Q_{2,0}(t)=\frac{105}{2}(t^2-3t^4)\\
Q_{2,1}(t)=\frac{105}{2}(\frac{1}{4}t^2-t^4+\frac{3}{4}t^6)\\
Q_{2,2}(t)=\frac{315}{64}(\frac{1}{6}t^2-\frac{1}{2}t^4+\frac{1}{2}t^6-\frac{1}{6}t^8).
\end{cases}
\end{equation}

\noindent Substituting the expressions from (\ref{77}) and (\ref{76}) into (\ref{9}) we obtain the approximation of $f$ for $n=2$. In Figure 3 we see the plot of the restriction of the approximation onto plane $y=3$.

\begin{figure}
\center
\includegraphics{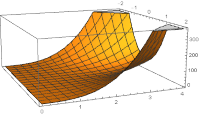}
\caption{the approximation of $f$ obtained by (\ref{9}) for $n=2$.}
\end{figure}

\textbf{Example 2}. Figure 4 provides an illustration of the iterative algorithm obtained by (\ref{9}) for a piecewise smooth function. We consider the following function

\begin{equation}\label{78}
f(x,y,z)=\begin{cases}
1 \,\,\,\,\,\, \texttt{for}\,\,\, x^2+y^2+(z-2)^2\leq1 \\
0\,\,\,\,\,\,  \texttt{for}\,\,\,
 x^2+y^2+(z-2)^2>1.\end{cases}
\end{equation}

\noindent It is obvious that on the section of the ball by plane $\{z=2\}$ we have function which takes $1$ in the unit disk with center $(0,0,2)$ and $0$ otherwise.

\noindent Now we are going to approximate the function using (\ref{9}) for $n=2$. For average $Mf(x,y,u)$ of $f$ over the sphere with center $(x,y,0)\in \{z=0\}$ and radius $u\geq 0$ we obtain
\begin{equation}\label{79}
Mf(x,y,u)=\begin{cases}
\frac{1}{2}-\frac{3+u^2+x^2+y^2}{4u\sqrt{4+x^2+y^2}} \,\,\,\texttt{for}\,\,\,\sqrt{4+x^2+y^2}-1 \leq u\leq \sqrt{4+x^2+y^2}+1\\
0\,\,\,\,\texttt{for}\,\,\,
 0<u\leq \sqrt{4+x^2+y^2}-1 \,\texttt{or}\, \sqrt{4+x^2+y^2}+1 \leq u  .\end{cases}
\end{equation}

\noindent Substituting the expressions from (\ref{77}) and (\ref{79}) into (\ref{9}) we obtain the approximation of $f$ for $n=2$. In Figure 4 we see the plot of the restriction of the approximation onto plane $z=2$.

\begin{figure}
\center
\includegraphics{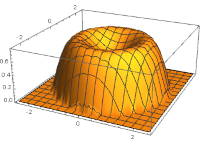}
\caption{}
\end{figure}

\end{document}